\documentstyle[amstex,amscd,12pt]{article}
\title{The Pontryagin rings of moduli spaces of
arbitrary rank holomorphic bundles over a Riemann surface}
\author{Richard Earl and Frances Kirwan}
\begin{document}
\setcounter{page}{1}
\newtheorem{prop}{PROPOSITION}
\newtheorem{lem}[prop]{LEMMA}
\newtheorem{cor}[prop]{COROLLARY}
\newtheorem{thm}[prop]{THEOREM}
\newtheorem{guess}{CONJECTURE}
\newtheorem{REM}[prop]{Remark}
\newenvironment{rem}{\begin{REM} \normalshape}{\end{REM}}
\newcommand{\Q}{\bold{Q}}
\newcommand{\mnd}{{\cal M}(n,d)}
\newcommand{\HS}{H^{*}}
\newcommand{\ar}{a_{r}}
\newcommand{\brk}{b_{r}^{k}}
\newcommand{\fr}{f_{r}}
\newcommand{\Z}{\bold{Z}}
\newcommand{\T}{\mbox{\bf \normalshape t}}
\newcommand{\D}{\mbox{\normalshape d}}
\newcommand{\Res}{\mbox{\normalshape Res}}
\newcommand{\ch}{\mbox{\normalshape ch}}
\newcommand{\mtc}{M_{\T}(c)}
\maketitle
\section{Introduction}
The cohomology of $\mnd$, the moduli space of stable holomorphic
bundles of coprime rank $n$ and degree $d$ and fixed determinant,
over a Riemann surface 
$\Sigma$ of genus $g \geq 2$, has been widely studied and from a wide
range of approaches. Narasimhan and Seshadri \cite{NS}
originally showed that the topology of $\mnd$ depends only on the
genus $g$ rather than the complex structure of $\Sigma$. An inductive
method to determine the Betti numbers of $\mnd$ was first given
by Harder and Narasimhan \cite{HN} and subsequently by Atiyah and
Bott \cite{AB}. The integral cohomology of $\mnd$
is known to have no torsion \cite{AB} and a set of
generators was found by Newstead
\cite{N} for $n=2$, and by Atiyah and Bott \cite{AB} for arbitrary
$n$. Much work and progress has been made recently in determining the
relations that hold amongst these generators, particularly in the rank
two, odd degree case which is now largely understood. A set of
relations due to Mumford in the rational cohomology ring of ${\cal
M}(2,1)$ is now known to be complete 
\cite{K}; recently several authors have
found a minimal complete set of relations for the `invariant' subring
of the rational cohomology of ${\cal M}(2,1)$ \cite{Z,B,KN,ST}.\\
\indent Unless otherwise stated all cohomology in this paper will have rational
coefficients.\\
\indent Let $V$ denote a normalised universal bundle over $\mnd \times
\Sigma$ \cite[p.582]{AB} and define classes
\begin{equation}
\ar \in H^{2r}(\mnd), \hspace{2mm} \brk \in H^{2r-1}(\mnd),\hspace{2mm}
\fr \in H^{2r-2}(\mnd),
\label{gen} 
\end{equation}
for $2 \leq r \leq n$ and $1 \leq k \leq 2g$ by setting
\[
c_{r}(V) = \ar \otimes 1 + \sum_{k=1}^{2g} \brk \otimes \alpha_{k} +
\fr \otimes \Omega \indent (2 \leq r \leq n)
\]
where $\alpha_{1},...,\alpha_{2g}$ is a fixed basis for
$H^{1}(\Sigma)$ and $\Omega$ is the standard generator for
$H^{2}(\Sigma)$. Atiyah and Bott \cite[Prop. 2.20]{AB} showed that the
rational cohomology ring $\HS(\mnd)$ is generated as a graded algebra
by the elements (\ref{gen}).\\
\indent The main results of this paper concern the vanishing of the
Pontryagin ring of $\mnd$ above a non-trivial degree.
\begin{thm}
\label{Pont}
The Pontryagin ring of $\mnd$ vanishes in degrees strictly greater
than $2n(n-1)(g-1).$
\end{thm} 
\indent The real dimension of $\mnd$ is $2(n^2-1)(g-1)$ and so Theorem
1 has consequence for $n \geq 2$ and $g \geq 2$. When $n =1$ or $g \leq
1$ the Pontryagin ring of $\mnd$ is trivial. 
\begin{thm}
\label{nonzero}
There exists a non-zero element of degree $2n(n-1)(g-1)$ in the
Pontryagin ring of $\mnd$.
\end{thm}
\indent When $n=2$, Theorem 1 is the first Newstead-Ramanan conjecture
\cite[p.344]{N}. In terms of the generators (\ref{gen}) above the
Newstead-Ramanan conjecture is equivalent to showing that
\[
(a_{2})^g = 0.
\]
This was first proved independently by Thaddeus \cite{T} and in
\cite[$\S$4]{K}. Subsequently it
has been proved in \cite{D, HS, JK, KN, W, WI2}. For the case $n=3$
Theorem 1 was also proved in \cite[$\S$5]{E}.\\
\indent Theorems 1 and 2 have recently been independently proved by
Jeffrey and Weitsman \cite{JW} for the arbitrary rank
case. When $n>2$ Theorems 1 and
2 are incompatible with a conjecture of Neeman \cite[p.458]{NE} which
stated that the Pontryagin ring of $\mnd$ should vanish in degree $2 g
n^2 - 4 g(n-1) +2$ and above.\\[\baselineskip]
\indent We shall prove these results by using formulas obtained in
\cite{JK,JK2} for the intersection pairings in $\HS(\mnd)$ between
cohomology classes represented as polynomials in the generators
$a_{r}, b_{r}^{k},f_{r}$. Knowing the intersection pairings of $\mnd$
one can of course (in
principal) determine the relations amongst the generators of
$\HS(\mnd)$, since by Poincar\'{e} duality an element $\zeta \in \HS(\mnd)$
of degree $p$ is zero if and only if
\[
\int_{\mnd} \eta \zeta = 0
\]
for every $\eta \in \HS(\mnd)$ of complementary degree
$2(n^{2}-1)(g-1) - p$. The results of \cite{JK,JK2} were inspired by
Witten's paper \cite{WI2} and use the principle of nonabelian
localization introduced in that paper and further developed in
\cite{GK,JK3,M}.\\[\baselineskip]
\indent The second Newstead-Ramanan conjecture states that the Chern
classes of ${\cal M}(2,1)$ also vanish above degree $4(g-1)$. This was
first proved geometrically by Gieseker \cite{G} and later by Zagier
\cite{Z} using Thaddeus' intersection pairings. In $\S$5 we give
explicit (though complicated) formulas for the pairings
\[
\int_{\mnd} \eta \cdot c(\mnd)(t)
\]
of arbitrary $\eta \in H^{*}(\mnd)$ with the Chern polynomial
$c(\mnd)(t)$ of $\mnd$. When $n=2$ and $d$ is odd, a proof of the
second Newstead-Ramanan conjecture may be easily rederived. 
Computer calculations for low values of $g$ and $n > 2$ suggest that
in general the Chern classes of $\mnd$ vanish above degree
$2n(n-1)(g-1)$.\\[\baselineskip]
\indent The last of the three Newstead-Ramanan conjectures states that
\[
\chi(\mnd, T\mnd) = 3 - 3 g.
\]
This was proved (for general $n$) by Narasimhan and Ramanan in
\cite{NR}. In fact they demonstrated the stronger result that
\[
H^{i}(\mnd,T\mnd) = \left\{ \begin{array}{ll} 3g-3 & i=1\\ 0 & i \neq 1.
\end{array} \right.
\]

\newpage
\section{Residue formulas for the intersection pairings in $\HS(\mnd)$}
\indent In \cite{JK2} formulas are given for the intersection pairings
in $\HS(\mnd)$ between cohomology classes expressed as polynomials
in the Atiyah-Bott generators (\ref{gen}). More precisely the evaluation $\int_{\mnd} \eta$
of the formal cohomology class
\begin{equation}
\eta = \exp (f_{2} + \delta_{3} f_{3} + \cdots + \delta_{n} f_{n} )
\prod_{r=2}^{n} \left( (a_{r})^{m_{r}} \prod_{k_{r}=1}^{2g}
(b_{r}^{k_{r}})^{p_{r,k_{r}}} \right) \label{eta}
\end{equation}
is considered where
\begin{itemize}
\item
$\delta_{3},...,\delta_{n}$ are formal nilpotent parameters,
\item
$m_{2},...,m_{n}$ are non-negative integers,
\item
$p_{r,k_{r}} \in \{0,1\}$ for $2 \leq r \leq n$ and $1
\leq k_{r} \leq 2g$.
\end{itemize}
Note that each $\brk$ has odd degree and hence $(\brk)^{2} = 0$. It
is sufficient to consider $\int_{\mnd} \eta$ for $\eta$ in the form
(\ref{eta}); by varying the integers $m_{r}$ and $p_{r,k_{r}}$ and
considering the coefficients of the monomials in $\delta_{3}, \ldots ,
\delta_{n}$ we may obtain the evaluation on the fundamental class
$[\mnd] \in H_{*}(\mnd)$ of any polynomial
in the generators (\ref{gen}), and hence the intersection pairing
\[
\langle \zeta, \xi \rangle = \int_{\mnd} \zeta \xi
\]
between any cohomology classes $\zeta,\xi \in \HS(\mnd)$ expressed
as polynomials in these generators.\\
\indent In \cite[Thm. 9.12]{JK2} the evaluation $\int_{\mnd} \eta$ of
$\eta$ on $[\mnd]$ is
equated to an iterated residue of a meromorphic function on the Lie
algebra
\[
\T = \{ \mbox{diag}(X_{1},...,X_{n}): X_{1} + \cdots +X_{n} = 0 \}
\]
of the standard maximal torus $T$ of $SU(n)$. There is a co-ordinate
system $(Y_{1},...,Y_{n-1})$ on $\T$ given by the simple roots $e_{1},
\ldots, e_{n-1}$ of $SU(n)$, i.e.
\[
Y_{j} = e_{j}(X) = X_{j} - X_{j+1} \mbox{ for } 1 \leq j \leq n-1,
\]
and the iterated residue is of the form
\[
\Res_{Y_{1}=0} \cdots \Res_{Y_{n-1}=0} \quad g (Y_{1}, \ldots Y_{n-1})
\]
where the variables $Y_{1},...,Y_{j-1}$ are held constant when
calculating $\Res_{Y_{j}=0}$ which is the usual residue at $0$ of a
meromorphic function in $Y_{j}$. If we use the Euclidean inner product
\[
\langle X,X \rangle = (X_{1})^{2} + \cdots + (X_{n})^{2}
\]
to identify $\T$ with its dual $\T^{*}$ then the simple roots
$e_{1},...,e_{n-1}$ correspond to generators
\[
\hat{e}_{j} = (0,...,0,1,-1,0,...,0)
\]
for the integer lattice of $\T$ (that is, the kernel of the
 exponential map from $\T$ to $T$). Let 
\[
\{\zeta_{j}^{k} : 1 \leq j \leq n-1, \quad 1 \leq k \leq 2g \}
\]
be the induced basis for $H^{1}(T^{2g})$.\\
\indent For $2 \leq r \leq n$ let $\sigma_{r}$ be the polynomial
function on $\T$ given by the $r$th elementary symmetric function in
$X_{1},...,X_{n}$ and for $X \in \T$ let
\[
q(X) = \sigma_{2}(X) + \delta_{3} \sigma_{3}(X) + \cdots + \delta_{n}
\sigma_{n}(X)
\]
where $\delta_{3},...,\delta_{n}$ are the formal nilpotent parameters
introduced above. We shall denote by $\D q_{X}:\T \to \bold{R}$ the
derivative of $q$ at $X \in \T$, so that
\[
\D q_{X} = \D(\sigma_{2})_{X} + \delta_{3} \D(\sigma_{3})_{X} + \cdots +
\delta_{n} \D(\sigma_{n})_{X}.
\]
The Hessian $\partial^{2}q_{X}$ of $q$ at $X$ is the symmetric
bilinear form on $\T$ given in any co-ordinate system by the matrix of
second partial derivatives of $q$ at $X$. Note that as $X_{1} + \cdots
+ X_{n} =0$ on $\T$ we have 
\begin{eqnarray*}
0 = (X_{1} + \cdots +X_{n})^{2} & = & \sum_{j=1}^{n} (X_{j})^{2} + 2
\sum_{i<j} X_{i}X_{j}\\
& = & \langle X, X \rangle + 2 \sigma_{2}(X)
\end{eqnarray*}
so that
\[
\sigma_{2}(X) = - \frac{1}{2} \langle X, X \rangle
\]
on $\T$. Then \cite[Thm. 9.12]{JK2} gives us the following formula for
$\int_{\mnd} \eta$.
\begin{thm}
\label{RES}
Let $\eta$ be the formal cohomology class given in (\ref{eta}). Then
$\int_{\mnd} \eta$ equals
\[
\frac{(-1)^{n(n-1)(g-1)/2}}{n!} \sum_{w \in W_{n-1}} \Res_{Y_{1}=0} \cdots \Res_{Y_{n-1}=0} \left[ \frac{ \exp \left\{
\D q_{X}(w \tilde{c}) \right\}
\prod_{r=2}^{n} \sigma_{r}(X)^{m_{r}} }{ {\cal D}(X)^{2g-2}
\prod_{j=1}^{n-1} \left( 1 - \exp (\D q_{X}(\hat{e}_{j})) \right)} \times
\right.
\]
\[
\times \left. \int_{T^{2g}} \exp \left\{ - 
\sum_{i,j=1}^{n-1} \sum_{k=1}^{g} \zeta_{i}^{k} \zeta_{j}^{k+g}
\partial^{2} q_{X}( \hat{e}_{i}, \hat{e}_{j}) \right\} \prod_{r=2}^{n}
\prod_{k_{r}=1}^{2g} \left( \sum_{j=1}^{n-1}
\D(\sigma_{r})_{X}(\hat{e}_{j}) \zeta_{j}^{k_{r}} \right)^{p_{r,k_{r}}}
\right]
\]
where
\[
{\cal D}(X) = \prod_{i<j}(X_{i}-X_{j})
\]
is the product of the
positive roots of $SU(n)$, where $\tilde{c} \in \T$ is the $n$-tuple
with $j$th entry
\[
\frac{d}{n} - \left[\frac{jd}{n} \right] + \left[\frac{(j-1)d}{n} \right],
\]
and where $W_{n-1} \cong S_{n-1}$ is the Weyl group of $SU(n-1)$
embedded in $SU(n)$ using the first $n-1$ co-ordinates.
\end{thm}
\begin{rem}
\label{newrem}
The formula of Theorem \ref{RES} is obtained by lifting the generators
$a_{r},b_{r}^{k}$ and $f_{r}$ of $\HS(\mnd)$ to $SU(n)$-equivariant
cohomology classes $\tilde{a}_{r},\tilde{b}^{k}_{r},\tilde{f}_{r}$ on
an `extended moduli space' (see \cite[$\S$4]{JK2}) with a Hamiltonian
$SU(n)$-action whose symplectic quotient is $\mnd$, and then using
`nonabelian localisation' \cite{JK2,WI} to localise to components of
the fixed point set of the maximal torus $T$. The restrictions of
$\tilde{a}_{r},\tilde{b}^{k}_{r},\tilde{f}_{r}$ to these components,
which are indexed by elements $\Lambda_{0}$ in the integer lattice of $\T$, and which can be identified with copies of $T^{2g}$, are
\[
\sigma_{r}(X), \indent \sum_{j=1}^{n-1} \D
(\sigma_{r})_{X}(\hat{e}_{j}) \zeta_{j}^{k}
\]
and
\[
\D (\sigma_{r})_{X}(\tilde{c} + \Lambda_{0}) + \sum_{i,j=1}^{n-1}
\sum_{k=1}^{g} \zeta_{i}^{k} \zeta_{j}^{k+g} \partial^{2}
(\sigma_{r})_{X}(\hat{e}_{i},\hat{e}_{j}) 
\]
respectively. As the extended moduli space is not compact, nonabelian
localisation cannot be applied directly, but instead one can exploit
the periodicity of the exponential map and the extended moduli space's
close relative 
\[
\mtc = \{ (h_{1},...,h_{2g},\Lambda) \in SU(n)^{2g} \times \T : 
\prod_{j=1} [h_{j},h_{j+g}] = c \exp(\Lambda) \}
\]
where 
\[
c = \exp(2 \pi i \tilde{c}) = e^{2 \pi i d / n} \mbox{diag}(1,1,\cdots
,1 ) \in SU(n).
\]
By \cite[Lemma 4.5]{JK2}, for any $\Lambda_{0}$ in the
integer lattice in $\T$ there
is a homeomorphism $s_{\Lambda_{0}}: \mtc \to \mtc$ defined by
\[
s_{\Lambda_{0}}(h,\Lambda) = (h, \Lambda + \Lambda_{0}).
\]
If as in \cite[$\S$ 3]{JK2} we represent
$T$-equivariant cohomology classes on a manifold $M$, acted on by the
torus $T$, by polynomial functions on the Lie algebra $\T$ of $T$ with
values in the De Rham complex of differential forms on $M$, then when
$X$ is the co-ordinate on $\T$, by (4.8) and (4.9) in \cite{JK2} we have
\begin{eqnarray*}
s^{*}_{\Lambda_{0}}(\tilde{a}_{r})(X) & = & \tilde{a}_{r}(X),\\
s^{*}_{\Lambda_{0}}(\tilde{b}_{r}^{k})(X) & = & \tilde{b}_{r}^{k}(X),
\end{eqnarray*}
and by \cite[Lemma 9.9]{JK2} 
\begin{equation}
\label{three}
s^{*}_{\Lambda_{0}} ( \tilde{f}_{2} + \sum_{r=3}^{n} \delta_{r}
\tilde{f}_{r} )(X) = (\tilde{f}_{2} + \sum_{r=3}^{n} \delta_{r}
\tilde{f}_{r} )(X) + \D \bar{q}_{X}(\Lambda_{0}) \label{cross}
\end{equation}
where
\[
\bar{q}(X) = \sigma_{2}(X) + \sum_{r=3}^{n} (-1)^r \delta_{r} \sigma_{r}(X).
\]
(Note that $\bar{q}$ was denoted by $q_{0}$ in \cite{JK2}.)\\
\indent Thus the result of applying $s^{*}_{\Lambda_{0}}$ to the
representative 
\[
\tilde{\eta} = \exp (\tilde{f}_{2} + \delta_{3} \tilde{f}_{3} + \cdots + \delta_{n} \tilde{f}_{n} )
\prod_{r=2}^{n} \left( (\tilde{a}_{r})^{m_{r}} \prod_{k_{r}=1}^{2g}
(\tilde{b}_{r}^{k_{r}})^{p_{r,k_{r}}} \right)
\]
of $\eta$ is 
\begin{equation}
\label{newchange}
s^{*}_{\Lambda_{0}}(\tilde{\eta}) = \tilde{\eta} \exp ( \D
\bar{q}_{X}(\Lambda_{0})).
\end{equation}
\indent Theorem \ref{RES} is now proved in \cite{JK2} using a version
of nonabelian localisation due to Martin \cite{M} and Guillemin and
Kalkman \cite{GK} which can be made to work in noncompact
settings. First one reduces to working with the symplectic
quotient of the extended moduli space by $T$ instead of by $SU(n)$ via
the arguments of \cite{M}. Then one compares the integral over the
$T$-quotient of the cohomology class induced by $\tilde{\eta}$ with
the integral of the class induced by $s^{*}_{\Lambda_{0}}$, for
$\Lambda_{0} = \hat{e}_{p}$ with $1 \leq p \leq n-1$, first using
(\ref{three}) and secondly using nonabelian localisation as in
\cite{GK} and \cite{M}.
\end{rem}
\begin{rem}
\label{epsilon}
The cohomology class $nf_{2}$ is represented by a symplectic form
$\omega$ on $\mnd$. If we replace $\omega$ by any non-zero scalar
multiple $\epsilon\omega$ then the proof of \cite[Thm. 9.12]{JK2} shows
that
\[
\int_{\mnd}\exp (\epsilon f_{2} + \delta_{3} f_{3} + \cdots + \delta_{n} f_{n} )
\prod_{r=2}^{n} \left( (a_{r})^{m_{r}} \prod_{k_{r}=1}^{2g}
(b_{r}^{k_{r}})^{p_{r,k_{r}}} \right)
\]
is given by the same iterated residue as in Theorem \ref{RES} above except that
$q(X)$ is replaced by
\[
q^{(\epsilon)}(X) = \epsilon\sigma_{2}(X) + \delta_{3} \sigma_{3}(X) +
\cdots + \delta_{n}\sigma_{n}(X)
\]
(see remarks 8.3 and 9.13 in \cite{JK2}). If we also multiply the
formal parameters $\delta_{3},...,\delta_{n}$ by $\epsilon$ then we
find that
\[
\int_{\mnd}\exp (\epsilon f_{2} + \epsilon\delta_{3} f_{3} + \cdots +
\epsilon \delta_{n} f_{n} ) \prod_{r=2}^{n} \left( (a_{r})^{m_{r}}
\prod_{k_{r}=1}^{2g} (b_{r}^{k_{r}})^{p_{r,k_{r}}} \right)
\]
is given by the same iterated residue as in Theorem \ref{RES} except that
$q(X)$ is now replaced by $\epsilon q(X)$.
\end{rem}

\newpage
\section{The Pontryagin Ring: Proof of Theorem \ref{Pont}}
Any symmetric polynomial in $X_{1},...,X_{n}$ can be expressed as a
polynomial in the elementary symmetric polynomials $\sigma_{1}
,...,\sigma_{n}$. Since $\sigma_{1}(X) = X_{1} + \cdots + X_{n}$ vanishes
on $\T$, any polynomial function on $\T$ which is symmetric in
$X_{1},...,X_{n}$ (or equivalently invariant under the action of the
Weyl group of $SU(n)$) can be expressed as a polynomial
$p(\sigma_{2}(X),\ldots,\sigma_{n}(X))$ and then represents the cohomology
class $p(a_{2},...,a_{n})$ on $\mnd$.
\begin{prop}
\label{Pontgen}
The Pontryagin ring of $\mnd$ is generated by the polynomials in
$a_{2},...,a_{n}$ represented by the elementary symmetric polynomials
in
\[
\{ (X_{i} - X_{j})^{2} : 1 \leq i < j \leq n \}.
\]
\end{prop}
{\bf Proof:} See \cite[Lemma 17]{E}.\\[\baselineskip]
Hence Theorem 1 is an immediate corollary of
\begin{thm}
\label{alt}
The subring of $\HS(\mnd)$ generated by $a_{2},...,a_{n}$ vanishes in
all degrees strictly greater than $2n(n-1)(g-1).$
\end{thm}
Since $\mnd$ is a compact manifold, by Poincar\'{e} duality
Theorem \ref{alt} is itself an immediate corollary of
\begin{prop}
\label{main}
Let $\eta$ be as given in (\ref{eta}). If $m_{2},...,m_{n}$ are
non-negative integers such that 
\[
\mbox{\normalshape deg} \prod_{r=2}^{n} (a_{r})^{m_{r}} = \sum_{r=2}^{n} 2 r m_{r}
> 2n(n-1)(g-1) 
\]
then $\int_{\mnd} \eta = 0$.
\end{prop}
{\bf Proof:} For $\epsilon \in {\bold R}$ let
\[
G(\epsilon) =  \int_{\mnd} \exp (\epsilon f_{2} + \epsilon \delta_{3}
f_{3} + \cdots + 
\epsilon \delta_{n} f_{n} ) 
\prod_{r=2}^{n} \left( (a_{r})^{m_{r}} \prod_{k_{r}=1}^{2g}
(b_{r}^{k_{r}})^{p_{r,k_{r}}} \right).
\]
We shall prove that $G(\epsilon) = 0$ for all $\epsilon \in {\bold R}$;
the result will then follow by taking $\epsilon =1$.\\
\indent First notice that 
\[
\mbox{deg }  \left( \prod_{r=2}^{n} \left( (a_{r})^{m_{r}} \prod_{k_{r}=1}^{2g}
(b_{r}^{k_{r}})^{p_{r,k_{r}}} \right) \right) = \sum_{r=2}^{n} 2 r
m_{r} + \sum_{r=2}^{n} \sum_{k_{r}=1}^{2g} (2 r -1) p_{r,k_{r}}.
\]
Also $f_{r}$ has degree $2r-2$, which is at least 2 for $2 \leq r \leq
n$, and the real dimension of $\mnd$ is $2(n^{2}-1)(g-1)$. Thus
$G(\epsilon)$ is a polynomial in $\epsilon$ of degree at most
\[
\frac{1}{2} \left( 2(n^{2}-1)(g-1) - \sum_{r=2}^{n} 2 r
m_{r} - \sum_{r=2}^{n} \sum_{k_{r}=1}^{2g} (2 r -1) p_{r,k_{r}}
\right).
\]
On the other hand Theorem \ref{RES} and Remark \ref{epsilon} show that
$G(\epsilon)$ is a non-zero $\epsilon$-independent constant multiple of 
\[
\sum_{w \in W_{n-1}} \Res_{Y_{1}=0} \cdots \Res_{Y_{n-1}=0} \left[
\frac{ \exp \left\{ \epsilon \D q_{X}(w\tilde{c}) \right\}
\prod_{r=2}^{n} \sigma_{r}(X)^{m_{r}} }{ {\cal D}(X)^{2g-2}
\prod_{j=1}^{n-1} \left(1 - \exp (\epsilon \D q_{X}(\hat{e}_{j})) \right)} \times
\right.
\]
\[
\times \left. \int_{T^{2g}} \exp \left\{ -\epsilon
\sum_{i,j=1}^{n-1} \sum_{k=1}^{g} \zeta_{i}^{k} \zeta_{j}^{k+g}
\partial^{2} q_{X}( \hat{e}_{i}, \hat{e}_{j}) \right\} \prod_{r=2}^{n}
\prod_{k_{r}=1}^{2g} \left( \sum_{j=1}^{n-1}
\D(\sigma_{r})_{X}(\hat{e}_{j}) \zeta_{j}^{k_{r}} \right)^{p_{r,k_{r}}}
\right].
\]
Now 
\[
q(X) = \sigma_{2}(X) + \delta_{3} \sigma_{3}(X) + \cdots + \delta_{n}
\sigma_{n}(X)
\]
on $\T$ where $\sigma_{2}(X) = - \langle X,X \rangle /2$ and
$\delta_{3},...,\delta_{n}$ are formal nilpotent parameters. Since
$\hat{e}_{j} \in \T$ corresponds, under the identification of $\T^{*}$
with $\T$ defined by the inner product, to the simple root $e_{j} \in
\T^{*}$ given by $e_{j}(X) = X_{j} - X_{j+1} = Y_{j},$ we have
$\D(\sigma_{2})_{X}(\hat{e}_{j}) = -Y_{j}$ and hence
\[
\D q_{X}(\hat{e}_{j}) = - Y_{j} + N_{j}
\]
where
\[
N_{j} = \delta_{3} \D(\sigma_{3})_{X}(\hat{e}_{j})+ \cdots +  \delta_{n}
\D(\sigma_{n})_{X}(\hat{e}_{j})
\]
is nilpotent. Moreover
\[
T(x) = \frac{x}{e^{x}-1}
\]
may be expressed as a formal power series in $x$ and then
\begin{eqnarray*}
\frac{1}{\exp(\epsilon \D q_{X}(\hat{e}_{j})) - 1} 
& = & \frac{1}{\exp(\epsilon N_{j} - \epsilon Y_{j}) -1}\\
& = & \frac{T(\epsilon N_{j} - \epsilon Y_{j})}{\epsilon N_{j} -
\epsilon Y_{j}}\\
& = &  - \frac{T(\epsilon N_{j} - \epsilon Y_{j})}{\epsilon Y_{j}}
\sum_{m=0}^{\infty} \left( \frac{N_{j}}{Y_{j}} \right)^{m} 
\end{eqnarray*}
The generators $\{ \zeta_{j}^{k} : 1 \leq j \leq n-1, 1 \leq k \leq 2g
\}$ for $H^{1}(T^{2g})$ are of degree 1 and hence anticommute, and all
have square 0. Their cup product is non-zero and hence spans the top
cohomology group $H^{2(n-1)g}(T^{2g})$. Thus as a function of
$\epsilon$ the integral
\[
\int_{T^{2g}} \exp \left\{ -\epsilon
\sum_{i,j=1}^{n-1} \sum_{k=1}^{g} \zeta_{i}^{k} \zeta_{j}^{k+g}
\partial^{2} q_{X}( \hat{e}_{i}, \hat{e}_{j}) \right\} \prod_{r=2}^{n}
\prod_{k_{r}=1}^{2g} \left( \sum_{j=1}^{n-1}
\D(\sigma_{r})_{X}(\hat{e}_{j}) \zeta_{j}^{k_{r}} \right)^{p_{r,k_{r}}}
\]
is a polynomial in $\epsilon$ which is divisible by $\epsilon$ to the
power
\[
\frac{1}{2} \left( \dim T^{2g} - \sum_{r=2}^{n} \sum_{k_{r}=1}^{2g}
p_{r,k_{r}} \right).
\]
Thus $G(\epsilon)$ is of the form 
\[
\Res_{Y_{1}=0} \cdots \Res_{Y_{n-1}=0} \left(
\frac{F(X,\epsilon)}{\cal{D}(X)^{2g-2} \prod_{j=1}^{n-1} (\epsilon
Y_{j})} \right)
\]
where, for $X \in \T$ and $\epsilon$ any real number, $F(X, \epsilon)$
is a formal power series in $\epsilon$ and a formal Laurent series in
the co-ordinates $Y_{1},...,Y_{n-1}$ on $\T$, and $F(X,\epsilon)$ is
divisible by $\epsilon$ raised to the power
\[
\frac{1}{2} \left( 2(n-1)g - \sum_{r=2}^{n} \sum_{k_{r}=1}^{2g}
p_{r,k_{r}} \right).
\]
Therefore $G(\epsilon)$ is a formal power series in $\epsilon$ which
is divisible by $\epsilon$ raised to the power
\[
\frac{1}{2} \left( 2(n-1)g - \sum_{r=2}^{n} \sum_{k_{r}=1}^{2g}
p_{r,k_{r}} \right) - (n-1) = 
\frac{1}{2} \left( 2(n-1)(g-1) - \sum_{r=2}^{n} \sum_{k_{r}=1}^{2g}
p_{r,k_{r}} \right)
\]
provided this is positive. On the other hand we saw earlier that
$G(\epsilon)$ is a polynomial in $\epsilon$ of degree at most
\[
\frac{1}{2} \left( 2(n^{2}-1)(g-1) - \sum_{r=2}^{n} 2 r
m_{r} - \sum_{r=2}^{n} \sum_{k_{r}=1}^{2g} (2 r -1) p_{r,k_{r}}
\right).
\]
Hence $G(\epsilon)$ must be identically zero unless
\[
2(n-1)(g-1) - \sum_{r=2}^{n} \sum_{k_{r}=1}^{2g}
p_{r,k_{r}} \leq  2(n^{2}-1)(g-1) - \sum_{r=2}^{n} 2 r
m_{r} - \sum_{r=2}^{n} \sum_{k_{r}=1}^{2g} (2 r -1) p_{r,k_{r}}.
\]
By hypothesis
\[
\sum_{r=2}^{n} 2 r m_{r} > 2 n (n-1) (g-1).
\]
So as $2r-1 \geq 1$ when $r \geq 1$ we have
\[ 
2(n-1)(g-1) - \sum_{r=2}^{n} \sum_{k_{r}=1}^{2g}
p_{r,k_{r}} >  2(n^{2}-1)(g-1) - \sum_{r=2}^{n} 2 r
m_{r} - \sum_{r=2}^{n} \sum_{k_{r}=1}^{2g} (2 r -1) p_{r,k_{r}}
\]
and hence $G(\epsilon)$ is identically zero. This completes the proof
of proposition \ref{main} and hence of Theorems \ref{Pont} and
\ref{alt}.

\section{The Pontryagin Ring: Proof of Theorem \ref{nonzero}}
By Proposition \ref{Pontgen} any symmetric polynomial in 
\[
\{(X_{i} - X_{j})^{2}: 1 \leq i < j \leq n\}
\]
represents an element of the Pontryagin ring of
$\mnd$. In particular the polynomial
\[
{\cal D}(X)^{2g-2} = \prod_{i<j}(X_{i}-X_{j})^{2g-2}
\]
represents an element $\eta_{0}$ of degree $2n(n-1)(g-1)$ in the
Pontryagin ring of $\mnd$. Thus Theorem \ref{nonzero} follows from the
following proposition.
\begin{prop}
\label{eta0}
If $\eta_{0} \in \HS(\mnd)$ is represented by ${\cal D}(X)^{2g-2}$
then
\[
\int_{\mnd} \eta_{0} \exp f_{2} \neq 0.
\]
\end{prop}
{\bf Proof:} By Theorem \ref{RES} $\int_{\mnd} \eta_{0} \exp f_{2}$ is a
non-zero constant multiple of 
\[
\sum_{w \in W_{n-1}}
\Res_{Y_{1}=0} \cdots \Res_{Y_{n-1}=0} \left[ \frac{ \exp \left\{
\D (\sigma_{2})_{X}(w\tilde{c}) \right\} }{
\prod_{j=1}^{n-1} \left( 1 - \exp (\D (\sigma_{2})_{X}(\hat{e}_{j})) \right)} \times
\right.
\]
\[
\times \left. \int_{T^{2g}} \exp \left\{ -
\sum_{i,j=1}^{n-1} \sum_{k=1}^{g} \zeta_{i}^{k} \zeta_{j}^{k+g}
\partial^{2} (\sigma_{2})_{X}( \hat{e}_{i}, \hat{e}_{j}) \right\}
\right]
\]
Now $\sigma_{2}$ is a quadratic form on $\T$ so $\D (\sigma_{2})_{X}$ is
linear in $X \in \T$ and $\partial^{2}(\sigma_{2})_{X}$ is independent
of $X$. Indeed we have already observed in the proof of Proposition \ref{main}
that $\sigma_{2}(X) = - \langle X,X \rangle /2$ so that
\[
\D(\sigma_{2})_{X}(\hat{e}_{j}) = - Y_{j}
\]
and 
\[
\D (\sigma_{2})_{X}(w\tilde{c}) = \beta_{1}(w)Y_{1} +
\beta_{2}(w) Y_{2} + \cdots + \beta_{n-1}(w) Y_{n-1}
\]
for constants $\beta_{j}(w)$, while
\[
\partial^{2}(\sigma_{2})_{X}(\hat{e}_{i},\hat{e}_{j}) = \left\{
\begin{array}{rl} -2 & \mbox{ if } i = j,\\
1 & \mbox{ if } i-j = \pm 1,\\
0 & \mbox{ if } |i-j|>1. \end{array} \right. 
\]
Since $\{ \hat{e}_{1},...,\hat{e}_{n-1} \}$ is a basis for $\T$ and
since $\HS(T^{2g})$ is a free exterior algebra on 
\[
\{\zeta_{j}^{k}: 1 \leq j \leq n-1, 1 \leq k \leq 2g\}
\]
it follows that
\[ 
\int_{T^{2g}} \exp \left\{ -
\sum_{i,j=1}^{n-1} \sum_{k=1}^{g} \zeta_{i}^{k} \zeta_{j}^{k+g} \right\}
\]
is a non-zero constant independent of $X \in \T$, which in fact equals
one. Hence
\[
\int_{T^{2g}} \exp \left\{ -
\sum_{i,j=1}^{n-1} \sum_{k=1}^{g} \zeta_{i}^{k} \zeta_{j}^{k+g}
\partial^{2} (\sigma_{2})_{X}( \hat{e}_{i}, \hat{e}_{j}) \right\}
\]
equals
\[
(\det \{ \partial^{2} (\sigma_{2})_{X}( \hat{e}_{i}, \hat{e}_{j})
\})^{g} 
\int_{T^{2g}} \exp \left\{ -
\sum_{i,j=1}^{n-1} \sum_{k=1}^{g} \zeta_{i}^{k} \zeta_{j}^{k+g} \right\}
= (-1)^{(n-1)g}n^{g}.
\]
\indent Thus $\int_{\mnd} \eta_{0} \exp f_{2}$ is
\[ (-1)^{(n-1)g}n^{g}
\sum_{w \in W_{n-1}}
\Res_{Y_{1}=0} \cdots \Res_{Y_{n-1}=0} \prod_{j=1}^{n-1} \left( \frac{
\exp\{\beta_{j}(w) Y_{j}\}}{1 - \exp\{-Y_{j}\} } \right) = (-1)^{(n-1)g}n^{g}
(n-1)!
\]
and the result follows. \indent $\Box$

\newpage
\section{Chern Classes}
The second Newstead-Ramanan conjecture states that the Chern classes
$c_{r}({\cal M}(2,1))$ vanish for $r > 2(g-1)$. This was
originally proved by Gieseker \cite{G} using a degeneration of the
moduli space and subsequently by Zagier \cite{Z} using Thaddeus'
pairings.\\

\indent The tangent bundle of $\mnd$ equals \cite[p.582]{AB}
\[
1 - \pi_{!}(\mbox{End} V)
\]
where $\pi: \mnd \times \Sigma \to \mnd$ is the first projection and
$V$ is a universal bundle over $\mnd \times \Sigma$. In
\cite[Prop. 10]{E} an expression for the Chern character of $V$ was
found in terms of the generators (\ref{gen}), and hence using the
Grothendieck-Riemann-Roch theorem an expression was determined for
$\ch(\mnd)$ in terms of these generators. If as in Remark \ref{newrem}
we lift these generators, and thus $\ch(\mnd)$, to $SU(n)$-equivariant
cohomology classes $\tilde{a}_{r},\tilde{b}_{r}^{k},\tilde{f}_{r}$ on
the extended moduli space and then restrict to the component indexed
by $\Lambda_{0} \in \T$ of the fixed point set of $T$, the result for
$\ch(\mnd)$ is  
\[
(1-g) + \sum_{i=1}^{n} \sum_{j=1}^{n} e^{X_{i}-X_{j}} (g-1 + \omega_{j} -
\omega_{i}) - \sum_{1 \leq i < j \leq n} \sum_{k=1}^{g}
Z_{i,j}^{k} Z_{i,j}^{k+g} (e^{X_{i}-X_{j}} + e^{X_{j}-X_{i}})
\]
where $\tilde{c} + \Lambda_{0} = (\omega_{1},...,\omega_{n}),$ and
\[
Z_{i,j}^{k} = - \zeta_{i-1}^{k} + \zeta_{i}^{k} + \zeta_{j-1}^{k} - \zeta_{j}^{k}, 
\]
with the understanding that $\zeta_{0}^{k} = \zeta_{n}^{k} = 0$. Hence
by \cite[Lemma 9]{E} the restriction of the lift of the Chern
polynomial $c(\mnd)(t) = \sum_{r \geq 0} c_{r}(\mnd) t^{r}$, which we
shall denote by $\tilde{c}(\mnd)(t)$, equals  
\[
\prod_{i<j} (1+(X_{i} - X_{j})t)^{g-1+\omega_{j}-\omega_{i}} (1+(X_{j} -
X_{i})t)^{g-1+\omega_{i}-\omega_{j}} \exp \left\{ \frac{-2t
\sum_{k=1}^{g} Z_{i,j}^{k} Z_{i,j}^{k+g}}{1 - (X_{i}-X_{j})^{2}
t^{2}} \right\}.
\]
\indent The formulas in Remark \ref{newrem} show that when
$\Lambda_{0} = \hat{e}_{p}$, 
\[
s_{\Lambda_{0}}^{*}(\tilde{c}(\mnd)(t)) = V_{p}(X,t) \cdot \tilde{c}(\mnd)(t)
\]
where
\[
V_{p}(X,t) = \prod_{q=1}^{n}
\frac{(1+(X_{q}-X_{p})t)(1+(X_{p+1}-X_{q})t)}{(1+(X_{p}-X_{q})t)(1+(X_{q}-X_{p+1})t)}.
\]
The proof of \cite[Thm. 9.12]{JK2} (see Remark \ref{newrem} above)
with this formula replacing \cite[Lemma 9.9]{JK2} (that is equations
(\ref{three}) and (\ref{newchange}) in Remark \ref{newrem}) shows that the
pairing $\int_{\mnd} \eta \cdot c(\mnd)(t)$ equals
$(-1)^{n(n-1)(g-1)/2}/n!$ times
\begin{eqnarray*} \sum_{w \in W_{n-1}}
\Res_{Y_{1}=0} \cdots \Res_{Y_{n-1}=0} 
\left\{
\frac{\exp \left\{ \D q_{X}(w\tilde{c})
\right\} \prod_{r=2}^{n} \sigma_{r}(X)^{m_{r}}}{
{\cal D}(X)^{2g-2} \prod_{p=1}^{n-1} \left( 1 - \exp ( \D
q_{X}(\hat{e}_{p})) V_{p}(X,t) \right)} \times
\right. \\
\nonumber \times  \left( \prod_{p \neq q}
(1+(X_{p}-X_{q})t)^{g-1+\tilde{c}_{q}-\tilde{c}_{p}} \right) \int_{T^{2g}}
\left [ 
\exp \left(
\sum_{i<j} \frac{ - 2 t  \sum_{k=1}^{g} Z_{i,j}^{k} Z_{i,j}^{k+g}
}{1 - (X_{i}-X_{j})^{2} t^{2}} \right) \right. \times \\
\nonumber \times \left. \left. \exp \left(
\sum_{i,j=1}^{n-1} \sum_{k=1}^{g} \zeta_{i}^{k}
\zeta_{j}^{k+g} \partial^{2} q_{X}(\hat{e}_{i},
\hat{e}_{j}) \right) \prod_{r=2}^{n} \prod_{k_{r}=1}^{2g} \left(
\sum_{j=1}^{n-1} \D (\sigma_{r})_{X}(\hat{e}_{j}) \zeta_{j}^{k_{r}}
\right)^{p_{r,k_{r}}} \right] \right\}.
\end{eqnarray*}
\newpage
When $n=2$ it readily follows that the resulting expression is a polynomial in
$t$ of degree at most $2g-2$ (see below). When $n \geq 3$, computer
calculations for low values of $n$ and $g$ suggest the following
(though we have thus far been unable to give proofs):
\begin{itemize}
\item $c_{r}(\mnd)=0$ for $r > n(n-1)(g-1)$.
\item $c_{n(n-1)(g-1)}(\mnd)$ is a non-zero multiple of $\eta_{0}$,
the cohomology class represented by ${\cal D}(X)^{2g-2}$ (see
Proposition \ref{eta0}).
\end{itemize}
For $n=2$ the results of Gieseker and Zagier may be duplicated as
follows. From \cite[p.144]{T} we need only consider those $\eta$ which
are invariant under the induced action of the mapping class group, so
for $\lambda \in {\bold C}$ let
\[
\eta = (a_{2})^{r} \exp \left\{ f_{2} + \lambda \sum_{k=1}^{g}
b_{2}^{k} b_{2}^{k+g} \right\}.
\]
For simplicity set $Y=Y_{1}$ and define
\[
F(k,s) = \Res_{Y=0} \left( \frac{[(1 - \lambda Y^{2}/2)(1-Y^{2}t^{2}) +
4 t]^{k}}{Y^{2s} (e^{Y/2} (1+Y t)^{2} - e^{-Y/2} (1-Y t)^2)} \right).
\]
Simplifying, we
see that $\int_{{\cal M}(2,1)} \eta c({\cal M}(2,1))(t)$ equals 
$(-1)^{g-1-r}2^{g-1-2r}F(g,g-1-r)$. We now claim that $F(k,s)$
is a rational function in $t$ of total degree at most $k+s-1$. We can show
by a simple induction that if $F(k,s)$ has total degree at most $k+s-1$ for
$k \leq K$ then $F(K+1,s)$ has total degree at most $K+s$. So consider 
$F(0,s)$. We may write  
\[
e^{Y/2} (1+Y t)^{2} - e^{-Y/2} (1-Y t)^2 = (1+ 4 t)Y \left(
1 + \sum_{i=1}^{\infty} r_{i}(t) Y^{2i} \right)
\]
where $r_{i}(t)$ is a rational function in $t$ of total degree 1. Hence
$F(0,s)$ is the coefficient of $Y^{2s}$ in 
\[
\frac{1}{(1+ 4 t)} \sum_{j=1}^{\infty}
\left( - \sum_{i=1}^{\infty} r_{i}(t) Y^{2i} \right)^{j}
\]
which is a rational function in $t$ of degree at most $s-1$, thus
proving the claim.\\
\indent Therefore $c({\cal M}(2,1))(t)$ is a polynomial in $t$ of
degree at most $2g-2$. It is in fact the case that the Chern
polynomial is of precisely this degree and we may find an expression
for $c_{2g-2}({\cal M}(2,1))$ as follows.\\

\indent Let $G(k,s)$ equal $F(k,s)$ modulo rational functions of
degree strictly less than $k+s-1$. Then
\[
G(0,s) = \frac{(-r_{1}(t))^{s}}{1+4t} = (-1)^{s} \frac{(t^2 +t/2 +
1/24)^{s}}{(1+4t)^{s+1}} = (-1)^{s} 2^{-2s-2} t^{s-1},
\]
and from the recurrence relation
\[
G(k+1,s) = 4 t G(k,s) - t^{2} G(k,s-1)
\]
we obtain
\[
G(k,s) = (-1)^{s} 2^{3k - 2s -2} t^{k+s-1}
\]
for $k \leq s.$ Using the above recurrence relation again we find
\[
G(s+1,s) = (-1)^{s} 2^{s} t^{2s} - t^{2} G(s,s-1) = (-1)^{s} t^{2s}
\sum_{i=0}^{s} 2^{i} = (-1)^{s} (2^{s+1}-1) t^{2s}.
\]
Hence the coefficient of $t^{2g-2}$ in $\int_{{\cal M}(2,1)} \eta
c({\cal M}(2,1))(t)$ equals $2^{g-1} (2^{g}-1)$.\\

\indent Let $\eta_{0}= (- 4 a_{2})^{g-1}$ be the class represented by
${\cal D}(X)^{2g-2} = Y^{2g-2}$. Then
\[
\int_{{\cal M}(2,1)} \eta_{0}\eta =\frac{(-1)^{g-1}}{2} \Res_{Y=0} \left[ \frac{
e^{-Y/2}}{1-e^{-Y}} \int_{T^{2g}} \exp \left\{- 2 \sum_{k=1}^{g} \zeta_{1}^{k}
\zeta_{1}^{k+g} \right\} \right] = (-2)^{g-1}.
\]
By Poincar\'{e} duality we find (cf. \cite[p.555]{Z})
\[
c_{2g-2}({\cal M}(2,1)) = (-1)^{g-1} (2^{g}-1) \eta_{0} = 2^{2g-2}(2^{g}-1) (a_{2})^{g-1}.
\]

Mathematical Institute, 24-29 St. Giles, Oxford OX1 3LB, England

\newpage
\begin{thebibliography}{WW}
\bibitem{AB} M.F.Atiyah and R.Bott {\em The Yang-Mills equations over
Riemann surfaces} Philos. Trans. Roy. Soc. London Ser. A {\bf 308}
(1982) 523-615.
\bibitem{B} V.Baranovsky {\em Cohomology ring of the moduli space of
stable vector bundles with odd determinant} Izv. Russ. Acad. Nauk.
{\bf 58 n4} (1994) 204-210.
\bibitem{D} S.K.Donaldson {\em Gluing techniques in the cohomology of
moduli spaces} in {\em Topological methods in modern mathematics}
(Proceedings of 1991 Stony Brook conference in honour of the sixtieth
birthday of J.Milnor) Publish or Perish.
\bibitem{E} R.A.Earl {\em The Mumford relations and the moduli of rank
three stable bundles} Compositio Math. (to appear)
\bibitem{G} D.Gieseker {\em A Degeneration of the Moduli Spaces of
Stable Bundles} J. Differential Geom. {\bf 19} (1984) 173-206.
\bibitem{GK} V.Guillemin and J.Kalkman {\em A new proof of the
Jeffrey-Kirwan localization theorem} J. Reine Agnew. Math. (to appear)
\bibitem{HN} G.Harder and M.S.Narasimhan {\em On the cohomology groups
of moduli spaces of vector bundles over curves} Math. Ann. {\bf 212}
(1975) 215-248. 
\bibitem{HS} R.Herrera and S.Salamon {\em Intersection numbers on
moduli spaces and symmetries of a Verlinde formula}
Comm. Math. Phys. (to appear)
\bibitem{JK3} L.C.Jeffrey and F.C.Kirwan {\em Localization for
nonabelian group actions} Topology {\bf 34} (1995) 291-327.
\bibitem{JK} L.C.Jeffrey and F.C.Kirwan {\em Intersection pairings in
moduli spaces of holomorphic bundles on a Riemann surface} Electronic
Research Announcements of the AMS 1 (1995) 57-71.
\bibitem{JK2} L.C.Jeffrey and F.C.Kirwan {\em Intersection theory on
moduli spaces of holomorphic bundles of arbitrary rank on a Riemann
surface} Ann. of Math. (to appear). 
\bibitem{JW} L.C.Jeffrey and J.Weitsman {\em Toric structures on the
moduli space of flat connections on a Riemann surface III: the higher
rank case} Duke Math. J. (to appear).
\bibitem{KN} A.D.King and P.E.Newstead {\em On the cohomology of the
moduli space of rank 2 vector bundles on a curve} Topology (to appear).
\bibitem{K} F.C.Kirwan {\em Cohomology rings of moduli spaces of
bundles over Riemann surfaces} J. Amer. Math. Soc. {\bf 5} (1992)
853-906.
\bibitem{M} S.K.Martin {\em Symplectic geometry and gauge theory}
Oxford D.Phil thesis (1997).
\bibitem{NR} M.S.Narasimhan and S.Ramanan {\em Deformations of the
moduli space of vector bundles over an algebraic curve} Ann. of
Math. {\bf 101} (1975) 391-417. 
\bibitem{NS} M.S.Narasimhan and C.S.Seshadri {\em Stable and unitary
vector bundles over an algebraic curve} Ann. of Math. {\bf 82} (1965) 540-567.
\bibitem{NE} A. Neeman {\em The topology of quotient varieties} Ann.
of Math. (2) {\bf 122} (1985) 419-459.
\bibitem{N} P.E.Newstead {\em Characteristic classes of stable
bundles of rank 2 over an algebraic curve} Trans. Amer. Math. Soc.
{\bf 169} (1972) 337-345.
\bibitem{ST} B.Siebert and G.Tian {\em Recursive relations for the
cohomology ring of moduli spaces of stable bundles} (Tr. J. of Math.)
{\bf 19} (1995) 131-144.
\bibitem{T} M.Thaddeus, {\em Conformal field theory and the cohomology
of the moduli space of stable  bundles} J. Differential Geom. {\bf
35} (1992) 131-149.
\bibitem{W} J.Weitsman {\em Geometry of the intersection ring of the moduli
space of
flat connections and the conjectures of Newstead and Witten} Topology
(to appear).
\bibitem{WI2} E.Witten {\em On quantum gauge theories in two
dimensions} Comm. Math. Phys. {\bf 141} (1991) 153-209.
\bibitem{WI} E.Witten {\em Two dimensional gauge theories revisited}
J. Geom. Phys. {\bf 9} (1992) 303-368.
\bibitem{Z} D.Zagier {\em On the cohomology of moduli spaces of rank 2
vector bundles over curves} (1995) Progress in Mathematics {\bf 129} {\em The
Moduli space of Curves} 533-563.
\end{thebibliography}
\end{document}